# Reduced Differential Transform Method for Solving Foam Drainage Equation(FDE)

Murat Gubes

Department of Mathematics, Karamanoğlu Mehmetbey University, Karaman/TURKEY

**Abstract:** Reduced Differental Transform Method (RDTM) which is one of the useful and effective numerical approximate method is applied to solve nonlinear time-dependent Foam Drainage Equation (FDE). Also, we compared the presented method with the famous Adomian Decomposition and Laplace Decomposition Technique. So, the results of giving examples demonstrated that RDTM is a powerful tool for solving nonlinear PDE, it can be applied very easily and it has less computational work than other existing methods like Adomian decomposition and Laplace decomposition.

## 1. Introduction

Most of the natural events, such as chemical, physical, biological, is modelled by a nonlinear equation. Besides exact solutions, we need its approximate solutions in terms of applicability. Therefore, a lot of approximate, numerical and analytic methods are developed and applied for nonlinear models [12-21].

Foam drainage which is a natural event described the process by which fluid flows out of a foam [2-3]. Since than many technological and industrial applications have been developed for foams, which include cleansing, water purification, minerals extraction [2-3].

More than ten years ago, given studies by Verbist and Weaire described the main features of both free drainage [5], [11], where liquid drains out of a foam due to gravity and forced drainage [11], where liquid is introduced to the top of a column of foam. In second state, a solitary wave of constant velocity is generated when liquid is added at a constant rate [10]. So forced foam drainage may be the best prototype for certain general phenomena described by nonlinear differential equations, particularly the type of solitary wave which is most familiar in tidal bores. The model developed by Verbist and Weaire idealizes the network of Plateau borders, through which the majority of liquid is assumed to drain, as a set of N identical pipes of cross section A, which is a function of position and time as in show [9], (1)-(2).

In this context, We investigate the famous time-dependent nonlinear forced channel-dominated foam drainage equation (FDE) to solve as numerical with reduced differential transform method (RDTM) [9]

$$\frac{\partial A(x,t)}{\partial t} + \frac{\partial}{\partial x}\left(A(x,t)^2 - \frac{\sqrt{A(x,t)}}{2}\frac{\partial A(x,t)}{\partial x}\right) = 0 \qquad (1)$$

where $x$ and $t$ are location and time respectively. In this case of forced drainage, it's solution is stated that [9]

$$A(x,t) = \begin{cases} c\tanh^2\left(\sqrt{c}(x-ct)\right), & x \leq ct \\ 0, & x \geq ct \end{cases} \qquad (2)$$

where $c$ is the velocity of wave front [11]. If we substitute $A(x,t) = u(x,t)^2$ and rearrange the equation (1) for basic form, then it can be written as follow

$$u(x,t)_t + 2u(x,t)^2 u(x,t)_x - u(x,t)_x^2 - \frac{1}{2}u(x,t)_{xx}u(x,t) = 0 \qquad (3)$$

with initial condition

$$u(x,0) = g(x). \qquad (4)$$

So, in this paper, we will consider and solve the equation (3) approximately and numerically for two different initial condition as (4).

## 2. Reduced Differential Transform Method (RDTM)

Let, suppose that $u(x,t)$ can be represented two variable function as a product of two single variable functions $f(x)\ and\ g(t)$ to show following manner

$$u(x,t) = f(x)g(t) \qquad (5)$$

From the similar meaning of definition of Differential Transform Method and its properties, we can write the transforming form of RDTM

$$u(x,t) = \sum_{i=0}^{\infty} F(i)x^i \sum_{j=0}^{\infty} G(j)t^j = \sum_{k=0}^{\infty} u_k(x)t^k \qquad (6)$$

where $u_k(x)$ is called $t$ dimensional spectrum function of $u(x,t)$.

If function $u(x,t)$ is analytic and differentiated continuously with respect to time $t$ and space $x$ in the domain of interest, then let

$$U_k(x) = \frac{1}{k!}\left[\frac{\partial^k}{\partial t^k}u(x,t)\right]_{t=0} \qquad (7)$$

Thus, from (7) it can be written the inverse transform of a sequence $\{U_k(x)\}_{k=0}^{\infty}$

$$u(x,t) = \sum_{k=0}^{\infty} U_k(x)t^k \qquad (8)$$

then combining (7) and (8), we obtain

$$u(x,t) = \sum_{k=0}^{\infty} \frac{1}{k!}\left[\frac{\partial^k}{\partial t^k}u(x,t)\right]_{t=0} t^k. \qquad (9)$$

If we consider the expressions (7),(8) and (9), it's clearly shown that the concept of the reduced differential transform is derived from the power series expansion. So, we give a table which included fundamental transformation properties of RDTM the following manner.

**Table 1: Basic transformations of RDTM for some functions.**

| Functional Form | Transformed Form |
|---|---|
| $u(x,t)$ | $U_k(x) = \frac{1}{k!}\left[\frac{\partial^k}{\partial t^k}u(x,t)\right]_{t=0}$ |
| $w(x,t) = u(x,t) \pm v(x,t)$ | $W_k(x) = U_k(x) \pm V_k(x)$ |

| $w(x,t) = \alpha u(x,t)$ | $W_k(x) = \alpha U_k(x)$, $\alpha$ constant |
|---|---|
| $w(x,t) = x^m t^n$ | $W_k(x) = x^m \delta(k-n)$, $\delta(k) = \begin{cases} 1, k = 0 \\ 0, k \neq 0 \end{cases}$ |
| $w(x,t) = x^m t^n u(x,t)$ | $W_k(x) = x^m U_{k-n}(x)$ |
| $w(x,t) = u(x,t)v(x,t)$ | $W_k(x) = \sum_{r=0}^{k} U_r(x) V_{k-r}(x) = \sum_{r=0}^{k} V_r(x) U_{k-r}(x)$ |
| $w(x,t) = \frac{\partial^r}{\partial t^r} u(x,t)$ | $W_k(x) = (k+1)\ldots(k+r) U_{k+r}(x)$ |
| $w(x,t) = \frac{\partial}{\partial x} u(x,t)$ | $W_k(x) = \frac{d}{dx} U_k(x)$ |
| $w(x,t) = \frac{\partial^2}{\partial x^2} u(x,t)$ | $W_k(x) = \frac{d^2}{dx^2} U_k(x)$ |

The proofs of Table 1 and the basic definitions of reduced differential transform method can be found in [23].

## 3. Implementations of RDTM and compared LDM-ADM

**3.1.** Consider the one-dimensional homogeneous forced foam drainage equation in (3) for the initial condition

$$u(x,0) = -\sqrt{c} \tanh(\sqrt{c} x). \tag{10}$$

Then, by using the basic properties of the reduced differential transformation, we can find the transformed form of equation (3) as

$$(k+1)U_{k+1}(x) = -2 \sum_{r=0}^{k} \sum_{s=0}^{k-r} U_r(x) U_s(x) \frac{d}{dx} U_{k-r-s}(x) + \sum_{r=0}^{k} \frac{d}{dx} U_r(x) \frac{d}{dx} U_{k-r}(x)$$
$$+ \frac{1}{2} \sum_{r=0}^{k} U_r(x) \frac{d}{dx^2} U_{k-r}(x) \tag{11}$$

and using the initial condition (4), we get

$$U_0(x) = -\sqrt{c}\tanh(\sqrt{c}x) \tag{12}$$

Now, put (12) into place (11), from hence we have the $U_k(x)$ values following

$$U_1(x) = \frac{c^2}{\cosh^2(\sqrt{c}x)}, U_2(x) = \frac{c^{\frac{7}{2}}\sinh(\sqrt{c}x)}{\cosh^3(\sqrt{c}x)}, U_3(x) = \frac{1}{3}\frac{c^5(2\cosh^2(\sqrt{c}x)-3)}{\cosh^4(\sqrt{c}x)}$$

$$U_4(x) = \frac{1}{3}\frac{c^{\frac{13}{2}}\sinh(\sqrt{c}x)(\cosh^2(\sqrt{c}x)-3)}{\cosh^5(\sqrt{c}x)}, \ldots \tag{13}$$

Thus, if we continue this process and also the inverse transformation of the set of $\{U_k(x)\}_{k=0}^{\infty}$ values gives approximate solution of (3) as

$$u(x,t) = \lim_{n\to\infty} \tilde{u}_n(x,t) = -\frac{1}{15}(15\sqrt{c}\sinh(\sqrt{c}x)\cosh^5(\sqrt{c}x) - 15c^2t\cosh^4(\sqrt{c}x)$$
$$-15c^{\frac{7}{2}}\sinh(\sqrt{c}x)t^2\cosh^3(\sqrt{c}x) + 15c^5t^3\cosh^2(\sqrt{c}x) - 10c^5t^3\cosh^4(\sqrt{c}x)$$
$$+15c^{\frac{13}{2}}t^4\sinh(\sqrt{c}x)\cosh(\sqrt{c}x) - 5c^{\frac{13}{2}}t^4\sinh(\sqrt{c}x)\cosh^3(\sqrt{c}x) - 15c^8t^5$$
$$+15c^8t^5\cosh^2(\sqrt{c}x) - 2c^8t^5\cosh^4(\sqrt{c}x))/\cosh^6(\sqrt{c}x) + \cdots \tag{14}$$

which converges the exact solution of (3) faster than adomian decomposition and laplace decomposition method.

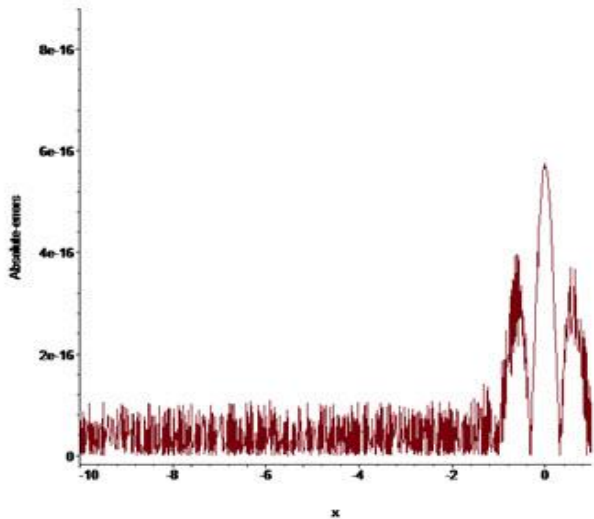

Figure 1: Absolute error of RDTM for $U_6$ at t=0.01, c=1 of 3.1

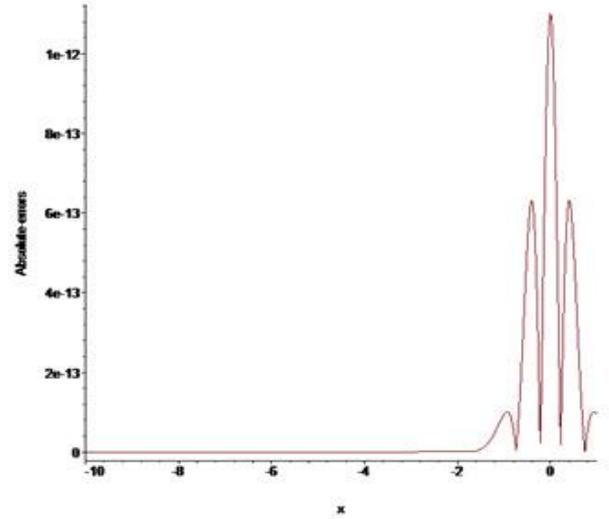

Figure1 2: Absolute error of RDTM for $U_6$ at t=0.01, c=2 of 3.1

**Table 2:** Comparison absolute errors of ADM, LDM and RDTM for ten terms approximation with $t = 0.1, c = 3$.

| x | ADM | LDM | RDTM |
|---|---|---|---|
| -10 | 0,202117817743509E-14 | 0,202117817743510E-14 | 0,202117817743509E-14 |
| -8 | 0,293257141972964E-14 | 0,293257141972664E-14 | 0,293257141972964E-14 |
| -6 | 0,45650705281285E-15 | 0,45650705281385E-15 | 0,45650705281285E-15 |
| -4 | 0,11116665360728E-12 | 0,11116665464828E-12 | 0,11116665360828E-12 |
| -2 | 0,851484003871439E-10 | 0,851483993870438E-10 | 0,851484003871439E-10 |
| 0 | 0,10317037658E-4 | 0,10317037658E-4 | 0,10317037658E-4 |

**Table 3:** Comparison absolute errors of ADM, LDM and RDTM for ten terms approximation with $t = 0.01, c = 3$.

| x | ADM | LDM | RDTM |
|---|---|---|---|
| -10 | 0,308694683847448E-15 | 0,308694683847448E-15 | 0,308694683847448E-15 |
| -8 | 0,507183534943173E-14 | 0,507183534943173E-14 | 0,507183534943173E-14 |
| -6 | 0,728297653568755E-15 | 0,728297653568755E-15 | 0,728297653568755E-15 |
| -4 | 0,695712017581145E-14 | 0,695712017481141E-14 | 0,695712017581145E-14 |
| -2 | 0,682573115077685E-15 | 0,682574115077685E-15 | 0,682573115077685E-15 |
| 0 | 0,1E-15 | 0,1E-15 | 0,1E-15 |

**Table 4:** Comparison absolute errors of ADM, LDM and RDTM for ten terms approximation with $t = 0.001, c = 3$.

| x | ADM | LDM | RDTM |
|---|---|---|---|
| -10 | 0,323329724865253E-16 | 0,323329724865253E-16 | 0,323329724865253E-16 |
| -8 | 0,300092135265223E-14 | 0,300092135265223E-14 | 0,300092135265223E-14 |
| -6 | 0,733106061807855E-14 | 0,733106061807855E-14 | 0,733106061807855E-14 |
| -4 | 0,483644638909495E-14 | 0,483644638909395E-14 | 0,483644638909495E-14 |
| -2 | 0,522768342580599E-14 | 0,522768342680599E-14 | 0,522768342580599E-14 |
| 0 | 0,4E-16 | 0,4E-16 | 0,4E-16 |

**Table 5:** The comparison of computation times which computed with Intel(R) Core (TM) i5-3230M CPU for 1.17 GHz between ADM, LDM and RDTM for equation 3.1.

| Iterations | Cpu times of ADM | Cpu times of LDM | Cpu times of RDTM |
|:---:|:---:|:---:|:---:|
| 5 steps | 0,747sec | 1,212sec | 0,488sec |
| 10 steps | 1,316sec | 2,813sec | 1,025sec |
| 15 steps | 2,67sec | 4,656sec | 1,846sec |
| 20 steps | 5,235sec | 9,478sec | 3,738sec |
| 25 steps | 9,755sec | 15,910se | 5,716sec |

**3.2.** Secondly, consider the nonlinear foam drainage equation (3) for different initial value such that

$$u(x,0) = -\frac{1}{2} + \frac{1}{1+e^x}. \tag{15}$$

Taking differential transform of (3) and the initial conditions (15) respectively, we obtain the same transformed form in (11) which has a different initial condition as follows

$$U_0(x) = -\frac{1}{2} + \frac{1}{1+e^x} \tag{16}$$

Then, substituting (16) into (11), we obtain the following $U_k(x)$ values,

$$U_1(x) = \frac{1}{4}\frac{e^x}{(1+e^x)^2}, U_2(x) = \frac{1}{32}\frac{e^x(-1+e^x)}{(1+e^x)^3}, U_3(x) = -\frac{1}{384}\frac{e^x(-1+4e^x-e^{2x})}{(1+e^x)^4}$$

$$U_4(x) = \frac{1}{6144}\frac{e^x(e^{3x}-1+11e^x-11e^{2x})}{(1+e^x)^5}, \tag{17}$$

$$U_5(x) = \frac{1}{122880}\frac{e^x(-66e^{2x}-1+26e^x+26e^{3x}-e^{4x})}{(1+e^x)^6}, \dots$$

Therefore, if we go on this process and the inverse transformation of the set of $\{U_k(x)\}_{k=0}^{\infty}$ values gives approximate solution of (3) as

$$u(x,t) = \lim_{n\to\infty} \tilde{u}_n(x,t) = \sum_{k=0}^{\infty} U_k(x)t^k \tag{18}$$

$$= -\frac{1}{122880} (61440 + 245760e^x + 307200e^{2x} - 245760e^{5x} - 61440e^{6x} + 7680e^{4x}t^2$$
$$+ 3840e^{5x}t^2 - 1920e^{3x}t^3 - 640e^{2x}t^3 - 640e^{4x}t^3 + 320e^{5x}t^3 + 20e^{5x}t^4$$
$$- 200e^{4x}t^4 + 200e^{2x}t^4 + 66e^{3x}t^5 - 26e^{2x}t^5 - 26e^{4x}t^5 + e^{5x}t^5$$
$$+ 122880e^{2x}t + 184320e^{3x}t + 122880e^{4x}t + 30720e^{5x}t - 7680e^{2x}t^2$$
$$+ e^x t^5 - 3840e^x t^2 + 320e^x t^3 - 20e^x t^4 + 30720e^x t - 307200e^{4x})$$
$$/(1+e^x)^6 + \cdots$$

which is similar solution with compared the adomian decomposition and laplace decomposition techniques and has also higher accuracy.

Table 6: Comparison between ADM, LDM and RDTM solutions of 3.2. at t=1.

| x   | ADM             | LDM             | RDTM            |
|-----|-----------------|-----------------|-----------------|
| -10 | 0,4999557229    | 0,4999557229    | 0,4999557229    |
| -9  | 0,4998796515    | 0,4998796515    | 0,4998796515    |
| -8  | 0,4996729266    | 0,4996729266    | 0,4996729266    |
| -7  | 0,4991114228    | 0,4991114228    | 0,4991114228    |
| -6  | 0,4975882790    | 0,4975882790    | 0,4975882790    |
| -5  | 0,4934713173    | 0,4934713173    | 0,4934713173    |
| -4  | 0,4824500775    | 0,4824500775    | 0,4824500775    |
| -3  | 0,4536908506    | 0,4536908506    | 0,4536908506    |
| -2  | 0,3833970310    | 0,3833970310    | 0,3833970310    |
| -1  | 0,2359453940    | 0,2359453940    | 0,2359453940    |
| 0   | 0,006249674499  | 0,006249674499  | 0,006249674499  |

Table 7: The comparison of computation times which computed with Intel(R) Core (TM) i5-3230M CPU for 1.17 GHz between ADM, LDM and RDTM for equation 3.2.

| Iterations | Cpu times of ADM | Cpu times of LDM | Cpu times of RDTM |
|------------|------------------|------------------|-------------------|
| 5 steps    | 0,745sec         | 1,182sec         | 0,491sec          |
| 10 steps   | 1,524sec         | 2,514sec         | 0,994sec          |
| 15 steps   | 3,292sec         | 4,499sec         | 2,245sec          |
| 20 steps   | 7,068sec         | 11,518sec        | 4,095sec          |
| 25 steps   | 14,153sec        | 23,32sec         | 7,681sec          |

# 4. Conclusions

Foam Drainage Equation is solved numerically by RDTM and compared with ADM-LDM. For example 3.1, given solution of RDTM is similar to ADM while it is better than LDM. At the same time given solution of RDTM provides the convergence more quickly than LDM and ADM as in show Table 5. Figure 1 and 2 indicate that our presented method is very effective and powerful.

Also, for example 3.2, presented method converges rapidly as well as obtained similar results with ADM and LDM as in Table 6 and 7.

# 5. References

[1] Weaire D, Hutzler S, Cox S, Alonso MD, Drenckhan D (2003) The fluid dynamics of foams. J Phys Condens Matter 15:65–72

[2] S. Hilgenfeldt, S.A.Koehler, H.A. Stone, Dynamics of coarsening foams: accelerated and self-limiting drainage, Phys.Rev.Lett, 20 (2001), 4704-4707.

[3] Koehler SA, Stone HA, Brenner MP, Eggers J (1998) Dynamics of foam drainage. Phys Rev E 58:2097–2106

[4] Weaire D, Hutzler S (2000) The physic of foams. Oxford University Press, Oxford

[5] Verbist G, Weaire D, Kraynik AM (1996) The foam drainage equation. J Phys Condens Matter 83:715–3731

[6] Helal MA, Mehanna MS (2007) The tanh method and Adomian decomposition method for solving the foam drainage equation. Appl Math Comput 190:599–609

[7] H. A. Stone, S. A. Koehler, S. Hilgenfeldt, and M. Durand, "Perspectives on foam drainage and the influence of interfacial rheology," *Journal of Physics Condensed Matter*, vol. 15, no. 1, pp. S283–S290, 2003.

[8] D.Weaire, S.Hutzler, G. Verbist, and E. A. J. Peters, "A review of foam drainage," *Advances in Chemical Physics*, vol. 102, pp. 315–374, 1997.

[9] D.Weaire, S. Findlay, and G. Verbist, "Measurement of foam drainage using AC conductivity," *Journal of Physics: Condensed Matter*, vol. 7, no. 16, pp. L217–L222, 1995.

[10] D. Weaire, S. Hutzler, N. Pittet, and D. Pardal, "Steady-state drainage of an aqueous foam," *Physical Review Letters*, vol. 71, no. 16, pp. 2670–2673, 1993.

[11] G. Verbist and D.Weaire, "Soluble model for foam drainage," *Europhysics Letters*, vol. 26, pp. 631–641, 1994.

[12] Yasir Khan, A method for solving nonlinear time-dependent drainage model, Neural Comput & Applic (2013) 23:411–415.

[13] Abdolhosein Fereidoon, Hessameddin Yaghoobi, and Mohammadreza Davoudabadi, Application of the Homotopy Perturbation Method for Solving the Foam Drainage Equation, Hindawi Publishing Corporation International Journal of Differential Equations Volume 2011, Article ID 864023, 13 pages doi:10.1155/2011/864023

[14] H. Jafari, V.D. Gejji. Adomian decomposition: a tool for solving a system of fractional differential equations. Mathematics Analysis Applications, 2005, 301: 508 - 518.

[15] H. Jafari, A. Yazdani, J. Vahidi, D.D. Ganji. Application of He's Variational Iteration Method for Solving Seventh Order Sawada-Kotera Equations. Applied Mathematical Sciences, 2008, 2: 471 - 477.

[16] A.M. Wazwaz. A study on linear and nonlinear Schrodinger equations by the variational iteration method. Chaos Solitons and Fractals, 2008, 37: 1136 - 1142.

[17] D. D. Ganji and M. Rafei, "Solitary wave solutions for a generalized Hirota-Satsuma coupled KdV equation by homotopy perturbation method," *Physics Letters. A*, vol. 356, no. 2, pp. 131–137, 2006.

[18] D. D. Ganji and A. Sadighi, "Application of He's homotopy-perturbation method to nonlinear coupled systems of reaction-diffusion equations," *International Journal of Nonlinear Sciences and Numerical Simulation*, vol. 7, no. 4, pp. 411–418, 2006.

[19] S. H. Mirmoradi, I. Hosseinpour, A. Barari and Abdoul R. Ghotbi, Analysis of Foam Drainage Problem Using Variational Iteration Method, Adv. Studies Theor. Phys., Vol. 3, 2009, no. 8, 283 - 292.

[20] A. Nikkar and M. Mighani, An Analytical Method to Analysis of Foam Drainage Problem, International Journal of Mathematical Sciences 2013 7(1).

[21] Majid Khan, Muhammad Asif Gondal, A new analytical solution of foam drainage equation by Laplace decomposition method, Journal of Advanced Research in Differential Equations, Vol. 2, Issue. 3, 2010, pp. 53-64.

[22] Y. Keskin, G. Oturanc, "Reduced Differential Transform Method for Partial Differential Equations", *International Journal of Nonlinear Sciences and Numerical Simulation*, 10(6) (2009) 741-749.

[23] Y. Keskin, Ph.D. Thesis, *Selcuk University* (to appear).

[24] Y. Keskin, G. Oturanc, "Reduced Differential Transform Method For Solving Linear And Nonlinear Wave Equations", *Iranian Journal of Science & Technology, Transaction A*, Vol. 34, No. A2, (2010).